\documentclass[12pt,reqno,oneside]{amsart}
 \usepackage{verbatim,amsmath,amssymb,cite,xspace}
\usepackage{color,cite,graphicx}

\theoremstyle{plain}
\newtheorem{theorem}{Theorem}[section]
\newtheorem{lemma}{Lemma}[section]

\theoremstyle{definition}
\newtheorem{definition}{Definition}[section]

\theoremstyle{remark}

\begin{document}
\title[]{Uniqueness of solutions to Navier Stokes equation with small initial data in $L^{3,\infty}(R^3)$}
\author{Hao Jia}
\maketitle
{\bf Abstract} In this short note we address a problem raised in \cite{Tsai}, concerning the uniqueness of solutions to Naiver Stokes equation with small initial data in $L^{3,\infty}(R^3)$, the Lorentz space. We prove uniqueness for such initial data.
\begin{section}{Introduction}
We consider the Navier Stokes equation in the whole space
\begin{eqnarray}\label{eq:NavierStokesEquation1}
\left.\begin{array}{rl}
\partial_tu-\Delta u+u\cdot\nabla u+\nabla p&=0\\
{\rm div~}u&=0\\
u(\cdot,0)&=u_0
\end{array}\right\}~~{\rm for~}(x,t)\in R^3\times(0,\infty).
\end{eqnarray}
If the initial data is small in the Lorentz space $L^{3,\infty}(R^3)$ (which is a critical space under the natural scaling of Navier Stokes equation), then by using perturbation method one can construct a global small solution $u$ with a number of regularity properties, such as
\begin{equation}
\sup_{t\ge 0}\left(t^{\frac{1}{4}}\|u(\cdot,t)\|_{L^6(R^3)}+t^{\frac{1}{8}}\|u(\cdot,t)\|_{L^4(R^3)}\right)<\infty.
\end{equation}
We must note that one can not construct even local in time solutions if the initial data is not small, in constrast to the case of initial data space $L^3(R^3)$ \footnote{If the initial data is scale-invariant, one does have existence of scale-invariant solutions, without any smallness condition on the initial data, see \cite{HaoSverak,Tsai}.}. Moreover, as pointed out in remarks below Lemma 4.1 of \cite{Tsai} , it seems that even when initial data is small in $L^{3,\infty}(R^3)$, the uniqueness of solutions has not been addressed before. There are of course many results concerning weak-strong uniqueness of solutions to Navier Stokes equation, see \cite{Leray,Prodi,Serrin2,LemRie} for example. For the initial data considered here, the situation is slightly more subtle. The reason roughly speaking is that $L^{3,\infty}(R^3)$ includes functions with a singularity of the order $O(\frac{1}{|x|})$, and singularity of this strength is precisely the borderline case of known regularity and uniqueness results for solutions to Navier Stokes equation, see \cite{HS} for more discussions. Thus previous results mostly implicitly assume $u_0$ has singularity of the form $o(\frac{1}{|x|})$ in appropriate senses. To illustrate the subtleties of $O(\frac{1}{|x|})$ type singularities further, we can consider the following semilinear heat equation in the whole space
\begin{equation}\label{eq:semilinearHeatEquation}
\partial_tu-\Delta u=u^3~~{\rm in~}R^3\times(0,\infty),
\end{equation}
with initial data $u_0(x)=\frac{\sigma}{|x|}$ and $\sigma>0$. We consider positive solution to equation (\ref{eq:semilinearHeatEquation}). This equation has similar scaling invariance as the Navier Stokes equations, it is invariant under the following scaling:
\begin{eqnarray}
u(x,t)&\rightarrow&u_{\lambda}(x,t):=\lambda u(\lambda x,\lambda^2 t),\\
u_0(x)&\rightarrow&u_{0\lambda}(x):=\lambda u_0(\lambda x).
\end{eqnarray}
It is proved in \cite{HarauxWeissler,Naito} that when $\sigma>0$ is small, there are two positive selfsimilar solutions $u(x,t)$ in the form $u(x,t)=\frac{1}{\sqrt{t}}\phi(\frac{x}{\sqrt{t}})$ with $\phi$ positive and radial. One of the solutions is small, and is given by perturbation method, the other solution is not small, and can not be detected by perturbation methods. If $\sigma>0$ is large, then there are no positive solutions to equation (\ref{eq:semilinearHeatEquation}). Notice that equation (\ref{eq:semilinearHeatEquation}) and equation (\ref{eq:NavierStokesEquation1}) have almost the same local existence theory based on perturbation methods. Here we show for Navier Stokes equation, such bizarre nonuniqueness can not happen. The difference is that Navier Stokes equation has energy inequality, which though simple, plays an important role in the properties of solutions. In our case the initial data space is $L^{3,\infty}(R^3)$ which has decay slower than the natural energy space $L^2(R^3)$, nonetheless we can still have a localized version of the energy inequality \cite{LemRie}. Such energy estimates give a priori bounds on weak solutions, and in particular imply all weak solutions are small if the initial data is small. This is our main observation to prove uniqueness. Our method is elementary and can be used to prove uniqueness with small initial data for more general initial data space $X$, which is scale invariant and locally stronger than $L^2_{loc}$. However we must note that our method can not be used to prove uniqueness for the space $BMO^{-1}(R^3)$ (the existence of small global solution has been proved in \cite{KoTa} with small initial data in $BMO^{-1}(R^3)$), since $BMO^{-1}(R^3)$ is not locally stronger than $L^2_{loc}$ (consequently weak solutions could still be large in energy space). Thus the following question is still open:

\medskip
\noindent
({\bf Q}) {\sl Suppose $u_0\in BMO^{-1}\cap L^2(R^3)$ with small $BMO^{-1}(R^3)$ norm, is the Leray solution with initial data $u_0$ unique?} 

\medskip
\noindent
Our paper is organized as follows. In section 2 we introduce the notion of Leray solution and provide some estimates when the initial data is small in $L^{3,\infty}(R^3)$; in section 3 we prove the uniqueness result.

\medskip
\noindent
{\bf Notations.} We use standard notations in the literature. For example, $B_R(x_0)$ denotes a ball with radius $R$ centered on $x_0$; $z_0=(x_0,t_0)$ is a point in the spacetime; $Q(R,z_0)=B_R(x_0)\times (t_0-R^2,t_0)$; for two vectors $u$ and $v$, $u\otimes v$ is a matrix with $(u\otimes v)_{ij}=u_iv_j$; for a matrix $a=(a)_{ij}$, ${\rm div~}a$ is a vector with $({\rm div~}a)_i=\partial_j a_{ij}$, where we assume Einstein's summation convention; we use $C$ to denote an absolute positive number, $C(\alpha,\beta,\dots)$ to denote a positive number depending on parameters $\alpha,\beta,\dots$; we adopt the convention that nonessential constants can vary from line to line; $u_0$ will be used as a divergence free vector field unless otherwise noted.
\end{section}

\begin{section}{Leray solutions and some estimates}
In \cite{LemRie} Lemarie-Rieusset introduced a very general notion of weak solution to the Navier Stokes equation. 
We recall the definition of Leray solutions in \cite{LemRie}, see also \cite{HaSv}.
\begin{definition}{(Leray solution)} A vector field
$u\in L^2_{loc}(R^3\times [0,\infty))$ is called a Leray solution to Navier-Stokes equations with initial data $u_0$ if it satisfies:

\smallskip
\noindent
i)   ${\rm ess}\sup_{0\leq t<R^2}\sup_{x_0\in R^3}\int_{B_R(x_0)}\frac{|u|^2}{2}(x,t)dx+\sup_{x_0\in R^3}\int_0^{R^2}\int_{B_R(x_0)}|\nabla u|^2dxdt<\infty$, and
$$\lim_{|x_0|\to \infty}\int_0^{R^2}\int_{B_R(x_0)}|u|^2(x,t)dxdt=0,$$
for any $R<\infty$.

\smallskip
\noindent
ii)  for some distribution $p$ in $R^3\times (0,\infty)$, $(u,p)$ verifies Navier Stokes equations
\begin{eqnarray}\label{eq:NavierStokesEquation}
 \left.\begin{array}{rl}
        \partial_tu-\Delta u+u\cdot\nabla u+\nabla p&=0\\
         \mbox{div~~}u&=0
       \end{array}\right\}\quad\mbox{in $R^3\times(0,\infty)$,}
\end{eqnarray}
in the sense of distributions and for any compact set $K\subseteq R^3$, $\lim_{t\to 0+}\|u(\cdot,t)-u_0\|_{L^2(K)}=0$.

\smallskip
\noindent
iii)  $u$ is suitable in the sense of Caffarelli-Kohn-Nirenberg, more precisely, the following local energy inequality holds:
\begin{equation}\label{eq:no2}
\int_0^{\infty}\int_{R^3}|\nabla u|^2\phi(x,t)dxdt\leq \int_0^{\infty}\int_{R^3}\frac{|u|^2}{2}(\partial_t\phi+\Delta \phi)+\frac{|u|^2}{2}u\cdot \nabla \phi+pu\cdot\nabla\phi dxdt
 \end{equation}
for any smooth $\phi\ge 0$ with ${\rm ~supp~~} \phi\Subset R^3\times (0,\infty)$. The set of all Leray solutions starting from $u_0$ will be denoted as $\mathcal{N}(u_0)$.
\end{definition}
We can calculate $p$ in the following way: $\forall (x,t)\in B_r(x_0)\times (0,t_{\ast})\subseteq R^3\times(0,\infty)$, take a smooth cutoff function $\phi$ with $\phi|_{B_{2r}(x_0)}=1$,  then there exists a function $p(t)$ depending only on $x_0,r,t,\phi$ (we suppress the dependence on $x_0,r,\phi$ in our notation) such that for $(x,t)\in B_r(x_0)\times(0,t_{\ast})$
\begin{equation}\label{eq:no3}
 p(x,t)=-\Delta^{-1}\mbox{div~div}(u\otimes u\phi)-\int_{R^3}\left(k(x-y)-k(x_0-y)\right)u\otimes u(y,t)\left(1-\phi(y)\right)dy+p(t)
\end{equation}
where $k(x)$ is the kernel of $\Delta^{-1}\mbox{div~div}$. \\
The right hand side is well defined since $u$ satisfies the estimates in i) and
\begin{equation}
|k(x-y)-k(x_0-y)|=O(\frac{1}{|x_0-y|^4}) {\rm ~~as~~} |y|\to \infty.
\end{equation}
The situation is similar to extending the domain of singular integrals to bounded functions, see for example \cite{LemRie} and \cite{Stein}.\\

For Leray solution $u\in \mathcal{N}(u_0)$, we have the following a priori estimates, first proved in \cite{LemRie}, see also a simpler proof in \cite{HaSv}. These estimates have played an important role in \cite{HaSv,RuSv}, see also \cite{Seregin}.\\
\begin{lemma}{\rm (A priori estimate for Leray solutions)}\label{lm:aprioriEstimatesForLeraySolution}\\
Let $\alpha=\sup_{x_0\in R^3}\int_{B_R(x_0)}\frac{|u_0|^2}{2}(x)dx<\infty$ for some $R>0$ and let $u$ be a Leray solution with initial data $u_0$. Then there exists some small absolute number $c>0$ such that for $\lambda$ satisfying $0<\lambda\leq c\min\{\alpha^{-2}R^2,1\}$, we have
\begin{equation}\label{eq:no4}
 {\rm ess}\sup_{0\leq t\leq \lambda R^2}\sup_{x_0\in R^3}\int_{B_R(x_0)}\frac{|u|^2}{2}(x,t)dx+\sup_{x_0\in R^3}\int_0^{\lambda R^2}\int_{B_R(x_0)}|\nabla u|^2(x,t)dxdt\leq C\alpha.
\end{equation}
\end{lemma}

\smallskip
\noindent
\textbf{Remarks:} Note that from the formula (\ref{eq:no3}) and the a priori estimate of $u$, we get the following estimate for $p$ which will be useful:
\begin{equation}
\sup_{x_0\in R^3}\int_0^{\lambda R^2}\int_{B_R(x_0)}|p-p(t)|^{\frac{3}{2}}dxdt\leq C \alpha^{\frac{3}{2}}R^{\frac{1}{2}}.
\end{equation}
Strictly speaking, one should really write $p_{x_0,R}(t)$ rather than just $p(t)$ in the last estimate. In other words, for a given $x_0, t$ and $R$ we need to choose an appropriate constant $p(t)=p_{x_0,R}(t)$ to satisfy the inequality. The reason is that the decay assumption on $u$ is too weak to imply decay of $p$, and the pressure may have a large mean value for large $|x_0|$.  However, the quantity entering the equation is $\nabla p$, which does not change in a given open set if we change $p$ in that set by a constant.  The estimates show that under our assumptions $p$ can be controlled up to such constants. The convention that $p(t)$ can depend on the corresponding $x_0$ and $R$ is used throughout the paper.\\

Using the above energy inequality, we can deduce the following estimate which will be useful below.
\begin{lemma}\label{lm:scaleinvariantestimate}
 Let $u$ be a Leray solution with initial data $u_0\in L^{3,\infty}(R^3)$. Let $p$ be the associated pressure. Then for $\forall r>0$,
\begin{eqnarray}
&&\int_0^{r^2}\int_{B_r(x_0)}|\nabla u|^2dxds+{\rm ess}\sup_{0\leq t\leq r^2}\int_{B_r(x_0)}\frac{|u|^2(x,t)}{2}dx\nonumber\\
&&\leq \frac{C\|u_0\|_{L^{3,\infty}(R^3)}^2r}{\sqrt{c\min\{\|u_0\|_{L^{3,\infty}(R^3)}^{-4},1\}}}{\rm ,~~and}\\
&&\int_0^{r^2}\int_{B_r(x_0)}|p-p(t)|^{\frac{3}{2}}dxds\leq \frac{C\|u_0\|_{L^{3,\infty}(R^3)}^3r^2}{c\min\{\|u_0\|_{L^{3,\infty}(R^3)}^{-4},1\}},
\end{eqnarray}
for any $x_0\in R^3$, where $c$ is the absolute small constant of Lemma \ref{lm:aprioriEstimatesForLeraySolution}.
\end{lemma}

\noindent
\textbf{Proof:} For each $r>0$, let $R=\frac{C_1r}{\sqrt{c\min\{\|u_0\|_{L^{3,\infty}(R^3)}^{-4},1\}}}>r$, with an absolute constant $C_1$ to be chosen below. We shall apply Lemma \ref{lm:aprioriEstimatesForLeraySolution} with this $R$. We have:
\begin{equation*}
\alpha=\sup_{x_0\in R^3}\int_{B_R(x_0)}|u_0|^2dx\leq C\|u_0\|_{L^{3,\infty}(R^3)}^2R.
\end{equation*}
We choose $C_1$ as the $C\ge 1$ in the above inequality, and set 
\begin{equation*}
\lambda=C^{-2}c\min\{\|u_0\|_{L^{3,\infty}(R^3)}^{-4},1\}\leq c\min\{\alpha^{-2}R^2,1\}. 
\end{equation*}
Note that by our choice of $R$ and $\lambda$, $\lambda R^2=r^2$, therefore from Lemma \ref{lm:aprioriEstimatesForLeraySolution}, the lemma follows.\\

The following version of $\epsilon$-regularity criteria of Caffarelli-Kohn-Nirenberg will be important for us in the sequel:
\begin{lemma}\label{lm:regularity}
Let $(u,p)$ be a suitable weak solution to NSE in $Q_R:=B_R(0)\times(-R^2,0)$ with $u\in L^{\infty}_tL^2_x(Q_R)\cap L^2_t\dot{H}^1(Q_R)$ and $p\in L^{\frac{3}{2}}(Q_R)$, in the sense that $(u,p)$ verifies NSE as distributions and they satisfy local energy inequality. Then there exists an absolute constant $\epsilon_0>0$, with the following property:\\
if $(R^{-2}\int_{Q_R}|u|^3dxdt)^{\frac{1}{3}}+(R^{-2}\int_{Q_R}|p|^{\frac{3}{2}}dxdt)^{\frac{2}{3}}\leq \epsilon_0$, then $\|\nabla^k u\|_{L^{\infty}(Q_{R/2})}\leq C_k R^{-k-1}$ for some constants $C_k$, $k=0,1,\dots$
\end{lemma}
A sketch of a short proof can be found for example in \cite{Lin}, a detailed one in \cite{LaSe}.\\

With the help of Lemma \ref{lm:scaleinvariantestimate} and \ref{lm:regularity}, we immediately get the following a priori estimate of any Leray solution $u$ to Navier Stokes equation with small initial data in $L^{3,\infty}(R^3)$.
\begin{lemma}\label{lm:smallenergysolution}
There exists a small positive number $\epsilon>0$ such that if $\|u_0\|_{L^{3,\infty}(R^3)}\leq \epsilon$, then for any Leray solution $u\in\mathcal{N}(u_0)$, we have
\begin{equation}
\sup_{t\ge 0} t^{\frac{1}{2}}\|u(\cdot,t)\|_{L^{\infty}(R^3)}\leq C(\epsilon)<\infty.
\end{equation}
Moreover, $C(\epsilon)\to 0+$ as $\epsilon\to 0+$.
\end{lemma}

\smallskip
\noindent
{\bf Proof and Remark.} For the proof, we only need to use the estimates from Lemma \ref{lm:scaleinvariantestimate} for $u$ in $B_R(x_0)\times (0,R^2)$ and apply Lemma \ref{lm:smallenergysolution} (provided $\epsilon$ is sufficiently small). We obtain
\begin{equation}
|u(x_0,R^2)|\leq C(\epsilon) R^{-1}.
\end{equation}
Take $R=\sqrt{t}$, we obtain the required estimate. 

\medskip
Define the uniformly locally bounded $L^p$ space $L^p_{uloc}$ ($p\ge 1$) as
\begin{equation}
L^p_{uloc}(R^3):=\{f\in L^p_{loc}:~\sup_{x_0\in R^3}\|f\|_{L^p(B_1(x_0))}<\infty\},
\end{equation}
with natural norm
\begin{equation}
\|f\|_{L^p_{uloc}(R^3)}:=\sup_{x_0\in R^3}\|f\|_{L^p(B_1(x_0))}.
\end{equation}
We will also need the following variation of the usual Young's inequality.
\begin{lemma}\label{lm:nonstandardconvolution}
Let $p\ge 1$, $p\leq q\leq \infty$, and let $\phi$ satisfy
\begin{equation}
\sup_{x\in R^3}(1+|x|)^4|\phi(x)|\leq M<\infty.
\end{equation}
Let $\phi_t(x)=t^{-\frac{3}{2}}\phi(\frac{x}{\sqrt{t}})$. Then for $0<t\leq 1$,
\begin{equation}
\|\phi_t\ast f\|_{L^q_{uloc}(R^3)}\leq C(p,q)Mt^{-\frac{3}{2}(\frac{1}{p}-\frac{1}{q})}\|f\|_{L^p_{uloc}(R^3)}.
\end{equation}
\end{lemma}

\smallskip
\noindent
{\bf Proof.} For any $x_0\in R^3$ we need to estimate $\|\phi_t\ast f\|_{L^q(B_1(x_0))}$. For any $f\in L^p_{uloc}$, write
\begin{equation}
f=f\chi_{B_2(x_0)}+f\chi_{B_2^c(x_0)}=f_1+f_2.
\end{equation}
We have by usual Young's inequality
\begin{equation}
\|\phi_t\ast f_1\|_{L^q(B_1(x_0))}\leq C(p,q)M t^{-\frac{3}{2}(\frac{1}{p}-\frac{1}{q})}\|f_1\|_{L^p(R^3)}\leq C(p,q)M t^{-\frac{3}{2}(\frac{1}{p}-\frac{1}{q})}\|f\|_{L^p_{uloc}(R^3)}.
\end{equation}
Note 
\begin{equation}
\|\phi_t\ast f_2\|_{L^q(B_1(x_0))}\leq C \|\phi_t\ast f_2\|_{L^{\infty}(B_1(x_0))}\leq CM \int_{B_2^c(x_0)}t^{-\frac{3}{2}}\frac{t^2}{|x_0-y|^4}f_2(y)dy\leq CM\|f\|_{L^p_{uloc}(R^3)}
\end{equation}
Thus the claim of the lemma follows.\\

Note from Lemma \ref{lm:smallenergysolution} and Lemma \ref{lm:scaleinvariantestimate} by interpolation we also have
\begin{equation}\label{eq:auxi1}
\sup_{x_0\in R^3,0\leq t\leq 1} t^{\frac{1}{4}}\|u(\cdot,t)\|_{L^4(B_1(x_0))}\leq C(\epsilon).
\end{equation}
However this estimate is not scale-invariant, we shall use the following improvement.
\begin{lemma}\label{lm:improvedregularity}
Let divergence free $u_0$ satisfy 
\begin{equation}
\|u_0\|_{L^{3,\infty}(R^3)}\leq \epsilon,
\end{equation}
for some sufficiently small positive number $\epsilon$. Then for any $u\in\mathcal{N}(u_0)$, 
\begin{equation}
\sup_{0\leq t\leq 1}t^{\frac{1}{8}}\|u(\cdot,t)\|_{L^4_{uloc}(R^3)}\leq C\epsilon.
\end{equation}
\end{lemma}

\smallskip
\noindent
{\bf Proof.} Denote for $0\leq t\leq 1$
\begin{equation}
\alpha(t):=\sup_{0\leq s\leq t}s^{\frac{1}{4}}\|u(\cdot,s)\|_{L^4_{uloc}(R^3)}.
\end{equation}
By (\ref{eq:auxi1}) we know that $\alpha(t)<\infty$ for all $1\ge t>0$. Note that $u$ satisfies the following integral equation
\begin{equation}
u(\cdot,t)=e^{\Delta t}u_0-\int_0^te^{\Delta (t-s)}P{\rm ~div~}(u\otimes u)(\cdot,s)ds,
\end{equation}
where $P$ is the Helmholtz projection operator to divergence free vector fields.\\
Since $\|u_0\|_{L^{3,\infty}(R^3)}\leq \epsilon$, we see
\begin{equation}
\sup_{t\ge 0}t^{\frac{1}{8}}\|e^{\Delta t}u_0\|_{L^4(R^3)}\leq C\epsilon.
\end{equation}
Note that the kernel $\phi_j(x)$ for $e^{\Delta}P\partial_j$ has the required property in Lemma \ref{lm:nonstandardconvolution}. Thus by using the estimate in Lemma \ref{lm:smallenergysolution} and the convolution inequality (with $p=q=4$) in Lemma \ref{lm:nonstandardconvolution}, we obtain
\begin{equation}
\alpha(t)\leq C\epsilon t^{\frac{1}{8}}+CC(\epsilon)\alpha(t).
\end{equation}
Let $\epsilon$ be sufficiently small such that $CC(\epsilon)<\frac{1}{2}$, we obtain
\begin{equation}
\alpha(t)\leq C\epsilon t^{\frac{1}{8}},
\end{equation}
which is exactly the claim in the lemma.
\end{section}

\begin{section}{Proof of main uniqueness result}
In this section we prove the following uniqueness result for Leray solutions with small initial data in $L^{3,\infty}(R^3)$.
\begin{theorem}
Suppose divergence free $u_0$ satisfy 
\begin{equation}
\|u_0\|_{L^{3,\infty}(R^3)}\leq \epsilon,
\end{equation}
for some sufficiently small positive number $\epsilon$. Then there is a unique Leray solution with initial data $u_0$.
\end{theorem}

\smallskip
\noindent
{\bf Proof.} The existence of Leray solution with initial data $u_0$ is well known, see for example \cite{LemRie}. Suppose there are two Leray solutions $u_1,~u_2\in \mathcal{N}(u_0)$. 
Taking $\epsilon$ sufficiently small, by the estimates in Lemma \ref{lm:improvedregularity}, we see
\begin{equation}
\sup_{1\ge t\ge 0}\left(t^{\frac{1}{8}}\|u_1(\cdot,t)\|_{L^4_{uloc}(R^3)}+t^{\frac{1}{8}}\|u_2(\cdot,t)\|_{L^4_{uloc}(R^3)}\right)\leq C\epsilon.
\end{equation}
Write
\begin{equation}
u_2=u_1+v.
\end{equation}
Then by substracting the equations for $u_1,~u_2$ we see $v$ satisfies the following integral equation
\begin{equation}
v(\cdot,t)=-\int_0^t e^{\Delta (t-s)}P{\rm ~div}(u_1\otimes v+v\otimes u_1+v\otimes v)(\cdot,s)ds.
\end{equation}
 Let
\begin{equation}
\alpha:=\sup_{0\leq t\leq 1}t^{\frac{1}{8}}\|v(\cdot,t)\|_{L^4(B_1(x_0))}.
\end{equation}
Then by the estimates of $u_1,~u_2$, Lemma \ref{lm:nonstandardconvolution} (with $p=2$ and $q=4$), we immediately obtain
\begin{equation}
\alpha\leq C\epsilon\alpha+C\alpha^2, ~~{\rm with~}C>1.
\end{equation}
Note that $\alpha\leq 2 C\epsilon$ by the estimates of $u_1,~u_2$. Take $\epsilon$ sufficiently small such that $C^2\epsilon\leq \frac{1}{8}$, we see from the above
that
\begin{equation}
\alpha\leq \frac{\alpha}{2},
\end{equation}
Thus $\alpha=0$ and $u_1\equiv u_2$ for $0\leq t\leq 1$. For $t\ge \frac{1}{2}$, $u_1$ and $u_2$ are bounded and have same value at $t=\frac{1}{2}$, thus $u_1\equiv u_2$ for all $t>0$ by the usual uniqueness results, see for example \cite{LemRie}. The theorem is proved.
\end{section}

\bigskip

{\bf Acknowledgement.} The author is grateful to V. Sverak for helpful comments and support.

\end{document}